\documentclass[12pt, letterpaper]{article}

\usepackage[utf8]{inputenc}
\usepackage[T1]{fontenc}      
\usepackage{amsmath, amssymb, amsthm}
\usepackage{geometry}
\usepackage{graphicx}
\usepackage{tikz}
\usepackage{url}
\usepackage{booktabs}
\usepackage{hyperref}

\theoremstyle{definition}
\newtheorem{definition}{Definition}[section]
\newtheorem{example}[definition]{Example}
\newtheorem{conjecture}{Conjecture}

\theoremstyle{plain}
\newtheorem{theorem}{Theorem}[section]

\title{Roller Coaster Permutations and Partition Numbers}
\author{William Adamczak \\
\small Email: wladamczak@gmail.com}
\date{\today}

\begin{document}

\maketitle

\begin{abstract}
This paper explores the partition properties of roller coaster permutations, a class of permutations characterized by maximizing the number of alternating runs in all subsequences. We establish a connection between the structure of these permutations and their partition numbers, defined as the minimum number of monotonic subsequences required to cover the permutation. Our main result provides a theoretical upper bound for the partition number of a roller coaster permutation of length $n$, given by $P_{max}(n) \le \lfloor\frac{\lceil\frac{n-2}{2}\rceil}{2}\rfloor + 2$. We further present experimental data for $n < 15$ that suggests this bound is nearly sharp.
\end{abstract}

\section{Introduction}

Roller coaster permutations, denoted as $RC_n$ for permutations of length $n$, were first introduced by Ahmed and Snevily \cite{Ahmed2013}. They are defined as permutations that maximize the total number of switches from ascending to descending (or vice versa) for the permutation itself and for all of its subpermutations simultaneously. More formally, these permutations achieve the greatest possible number of local extrema (peaks and valleys) across all possible subsequences.

It is important to note that the term "rollercoaster" has recently been used by Biniaz et al. \cite{Biniaz2019} to describe sequences where every run has a length of at least 3. However, throughout this paper, we strictly adhere to the original definition by Ahmed and Snevily regarding maximizing alternating subsequences.

Recent work by Netto \cite{Netto2022} and Botler \& Netto \cite{Botler2021} has further explored these permutations using Integer Linear Programming to find new bounds. In his work, Netto proposed a counterexample to the proof of the alternating property established in previous work by the author \cite{Adamczak2017}. Specifically, for the permutation $\pi=7,3,5,6,4,8,1$, Netto suggests that swapping 5 and 6 (yielding $\pi^{\prime}=7,3,6,5,4,8,1$) decreases the alternation value, appearing to contradict the claim that swapping adjacent elements in a monotonic triple increases the value. However, the theorem in \cite{Adamczak2017} prescribes swapping the \emph{first two} elements of the monotonic triple (in this case, the 3 and the 5), not the second and third. As the counterexample relies on a swap operation distinct from the one specified in the theorem, the original proof strategy remains valid.

Roller coaster permutations share deep connections with the theory of pattern avoidance. Stankova \cite{Stankova1994} classified permutations based on forbidden subsequences and established that specific classes of these pattern-avoiding permutations are subsets of $RC_n$. Consequently, the properties of roller coaster permutations are intrinsically linked to stack-sortable permutations, as seen in Egge \& Mansour \cite{Egge2004}, where the connection between forbidden subsequences and stack sortability is established.

The primary focus of this paper is the partition number of these permutations. The partition number corresponds to the minimum number of monotonic subsequences required to partition $\pi$. This metric has significant applications in computer science, specifically corresponding to the number of piles formed in the card game Patience Sorting \cite{Aldous1999}. While calculating the partition number for increasing subsequences alone is efficiently solvable via the Robinson-Schensted correspondence, Wagner \cite{Wagner1984} showed that the problem of partitioning into sets of \emph{both} increasing and decreasing subsequences is generally NP-hard. This computational difficulty underscores the value of establishing theoretical bounds for specific classes of permutations, such as roller coasters.

While Kezdy, Snevily \& Wang \cite{Kezdy1996} explored partitions by associating a graph to a permutation, our approach relies on the underlying alternating structure and the strict relative positions of entries forced by the definition of roller coaster permutations.

\section{Background}

To analyze the partition properties of roller coasters, we first establish the necessary terminology regarding permutation patterns and structure.

\begin{definition}
A \textbf{permutation} of length $n$ is an ordered rearrangement on the set $\{1,2,3,\dots,n\}$. The collection of all such permutations is denoted $S_{n}$.
\end{definition}

The central object of this study is the roller coaster permutation, originally described by Ahmed and Snevily \cite{Ahmed2013} to capture the concept of maximality in direction changes.

\begin{definition}
A \textbf{Roller Coaster permutation} is a permutation that maximizes the number of changes from increasing to decreasing over itself and all of its subsequences, simultaneously. The set of all such permutations in $S_n$ is denoted $RC_n$.
\end{definition}

For small values of $n$, the sets $RC_n$ have been explicitly enumerated in \cite{Ahmed2013}. It is worth noting that $|RC_n|$ is generally larger for odd $n$ than for even $n$, a property likely arising from the symmetric flexibility of the central element in odd-length permutations.

The structure of these permutations relies heavily on the concept of alternation.

\begin{definition}
An \textbf{alternating permutation} is a permutation $\pi$ such that $\pi_{1} < \pi_{2} > \pi_{3} \dots$. Conversely, a \textbf{reverse alternating permutation} satisfies $\pi_{1} > \pi_{2} < \pi_{3} \dots$.
\end{definition}

\begin{example}
The permutation $\pi = \{4,3,7,1,8,2,6,5\}$ is reverse alternating. The first entry is greater than the second, and this alternating pattern continues throughout the sequence.
\end{example}

Prior work by the author \cite{Adamczak2017} and Ahmed \& Snevily \cite{Ahmed2013} established several structural theorems that are critical for deriving our partition bound. Most notably, roller coaster permutations are closed under basic symmetry operations and adhere to strict alternating patterns.

\begin{theorem}\label{thm:ends}
For any $\pi \in RC_{n}$:
\begin{enumerate}
    \item The reverse and complement of $\pi$ are also in $RC_{n}$ \cite{Ahmed2013}.
    \item $\pi$ is strictly alternating or reverse alternating \cite{Adamczak2017}.
    \item The endpoints satisfy $|\pi_{1}-\pi_{n}|=1$ \cite{Adamczak2017}.
\end{enumerate}
\end{theorem}

These structural constraints enforce a strict separation of values based on index parity. This separation is the key feature that allows us to partition the permutation efficiently.

\begin{theorem} \label{thm:separation}
For $\pi \in RC_{n}$:
\begin{enumerate}
    \item If $\pi$ is alternating, then $\pi_{i} > \max(\pi_{1}, \pi_{n})$ for all even $i$.
    \item If $\pi$ is reverse alternating, then $\pi_{i} > \max(\pi_{1}, \pi_{n})$ for all odd $i$ (excluding endpoints).
\end{enumerate}
\end{theorem}

While we have established that roller coaster permutations are alternating, they satisfy a much stronger structural condition known as being \emph{recursively alternating}, as defined in \cite{Adamczak2017}.

\begin{definition}
Let $\pi[2, n-1]$ denote the restriction of $\pi$ to the internal indices. Let $\pi_{odd}$ be the restriction of $\pi[2, n-1]$ to odd indices, and $\pi_{even}$ be the restriction to even indices. A permutation $\pi$ is said to be \textbf{recursively alternating} if:
\begin{enumerate}
    \item $\pi$ is alternating or reverse alternating.
    \item $\pi_{odd}$ and $\pi_{even}$ are themselves alternating or reverse alternating.
    \item This property holds recursively for deeper restrictions (i.e., for the odd/even restrictions of $\pi_{odd}$, etc.).
\end{enumerate}
Equivalently, the restriction of $\pi[2, n-1]$ to indices congruent modulo $2^k$ is alternating for any $k$.
\end{definition}

\begin{theorem} \label{thm:recursive}
If $\pi \in RC(n)$, then $\pi$ is recursively alternating \cite{Adamczak2017}.
\end{theorem}

This theorem is powerful because it guarantees that the subsequences formed by taking every second element (or every fourth, etc.) retain the alternating structure.

To illustrate the structure of roller coaster permutations consider the following example:

\begin{example}
Consider the permutation $\pi = \{5,3,7,1,8,2,6,4\}$. A graphical representation is provided in Figure \ref{fig:value_separation}.
Observe the following structural properties:
\begin{itemize}
    \item The endpoints are $5$ and $4$, satisfying the condition $|\pi_1 - \pi_n| = 1$ established in Theorem \ref{thm:ends}.
    \item The internal elements are strictly separated by the endpoints. The elements $\{7, 8, 6\}$ form a "Top Set" where every value is greater than the endpoints. Conversely, the elements $\{3, 1, 2\}$ form a "Bottom Set" where every value is less than the endpoints.
    \item This corresponds precisely to the index parity: the values in the Top Set occur at odd indices ($3, 5, 7$), while the values in the Bottom Set occur at even indices ($2, 4, 6$).
    \item The alternating (or reverse-alternating) behavior of the Top and Bottom sets is also readily seen in this diagram.
\end{itemize}

\begin{figure}[h]
    \centering
    \begin{tikzpicture}[scale=1.0]
        \tikzstyle{vertex}=[circle, draw=black, fill=white, inner sep=2pt]
        \tikzstyle{permline}=[black, thin]
        \tikzstyle{guide}=[gray, dashed, thin]
        
        \coordinate (p1) at (1,5);
        \coordinate (p2) at (2,3);
        \coordinate (p3) at (3,7);
        \coordinate (p4) at (4,1);
        \coordinate (p5) at (5,8);
        \coordinate (p6) at (6,2);
        \coordinate (p7) at (7,6);
        \coordinate (p8) at (8,4);

        \draw[guide] (0.5, 5) -- (8.5, 5) node[right, black] {End $\pi_1=5$};
        \draw[guide] (0.5, 4) -- (8.5, 4) node[right, black] {End $\pi_n=4$};

        \draw[permline] (p1)--(p2)--(p3)--(p4)--(p5)--(p6)--(p7)--(p8);

        \draw[thick, red] plot [smooth cycle, tension=0.7] coordinates {
            (2.6, 6.8) (3, 7.5) (5, 8.5) (7, 6.5) (7.4, 5.8) (5, 6.5) 
        };
        \node[red] at (5, 7.5) {Top ($O$)};

        \draw[thick, blue] plot [smooth cycle, tension=0.7] coordinates {
            (1.6, 3.2) (2, 2.5) (4, 0.5) (6, 1.5) (6.4, 2.2) (4, 2.5)
        };
        \node[blue] at (4, 1.5) {Bottom ($E$)};

        \foreach \i in {1,...,8} {
            \node[vertex] at (p\i) {};
        }
        
        \node[above] at (p1) {5};
        \node[below] at (p2) {3};
        \node[above] at (p3) {7};
        \node[below] at (p4) {1};
        \node[above] at (p5) {8};
        \node[below] at (p6) {2};
        \node[above] at (p7) {6};
        \node[below] at (p8) {4};

    \end{tikzpicture}
    \caption{Graphical representation of $\{5,3,7,1,8,2,6,4\}$. The horizontal lines at $y=4$ and $y=5$ demonstrate the strict value separation. Elements at even indices (Bottom Set) are strictly below the endpoints, while elements at odd indices (Top Set) are strictly above.}
    \label{fig:value_separation}
\end{figure}
\end{example}

Finally, we define the metric of interest for this paper.

\begin{definition}
A \textbf{partition} of a permutation is a set of disjoint monotonic subsequences (increasing or decreasing) that cover the permutation. The \textbf{partition number}, denoted $P(\pi)$, is the minimum size of such a set.
\end{definition}

We denote $P_{max}(n)$ as the maximum partition number among all $\pi \in RC_{n}$. Our goal is to bound this value.

\section{Results}

\begin{theorem}
For any $\pi \in RC_{n}$, the partition number $P(\pi)$ is bounded above by:
\[ P_{max}(n) \le \lfloor\frac{\lceil\frac{n-2}{2}\rceil}{2}\rfloor + 2 \]
\end{theorem}

\begin{proof}
Without loss of generality, assume $\pi \in RC_{n}$ is reverse alternating (i.e., $\pi_{1} > \pi_{2} < \pi_{3} \dots$). If $\pi$ were alternating, we could consider its complement, which is also in $RC_{n}$ (Theorem 2.9) and is reverse alternating. Since the partition number is invariant under complementation, the bound derived here applies to all roller coaster permutations.

Let $\pi_{1}$ and $\pi_{n}$ be the endpoints. We separate the internal elements into two sets based on index parity:
\[ E = \{ \pi_i \mid 2 \le i \le n - 1, i \text{ is even} \} \quad (\text{"Bottom" set}) \]
\[ O = \{ \pi_i \mid 2 \le i \le n - 1, i \text{ is odd} \} \quad (\text{"Top" set}) \]

By Theorem \ref{thm:separation}, strict value separation occurs: for all $e \in E$ and $o \in O$, we have $\pi_{e} < \min(\pi_{1}, \pi_{n}) < \max(\pi_{1}, \pi_{n}) < \pi_{o}$. This separation allows us to partition $E$ and $O$ separately and then stitch them together.

First, we decompose $E$ and $O$ into the minimum number of contiguous increasing runs. By Theorem \ref{thm:recursive}, $\pi$ is recursively alternating; thus, the subsequences defined by $E$ and $O$ are themselves alternating (or reverse alternating). Noting that $E$ has $\lceil\frac{n-2}{2}\rceil$ elements and the number of runs can exceed half of this by at most 1 due to the alternating behavior, the number of runs required for $E$ is bounded by:
\[ \text{Runs}(E) \le \lfloor\frac{\lceil\frac{n-2}{2}\rceil}{2}\rfloor + 1 \]

We now construct the partitions by pairing these runs. Let $R_{top}^{(k)}$ be the $k$-th run in $O$ and $R_{bottom}^{(k)}$ be the $k$-th run in $E$.
\begin{enumerate}
    \item \textbf{Index Constraints:} When partitioning a forward or reverse alternating permutation into contiguous increasing runs, the $k$-th run has an earliest start at index $2k-2$ (for $k>1$) and a latest finish at index $2k$. Due to the reverse alternating structure of $\pi$, the indices of the bottom set ($E$) precede the corresponding indices of the top set ($O$). Specifically, the $k$-th index of the bottom partition generally precedes the $k+1$-th index of the top partition.
    
    \item \textbf{Pairing Strategy:}
    \begin{itemize}
        \item The \textbf{Start Point} $\pi_1$ pairs with the first run on the top, $R_{top}^{(1)}$. This is valid because $\pi_1$ is at index 1 (preceding $O$'s indices) and $\pi_1 < O$.
        \item The $k$-th run on the bottom, $R_{bottom}^{(k)}$, pairs with the $(k+1)$-th run on the top, $R_{top}^{(k+1)}$. The value separation condition ($E < O$) ensures validity, and the index logic described above ensures the bottom run finishes before the top run starts.
    \end{itemize}

    \item \textbf{Endpoints:}
    If the number of runs in $E$ exceeds the number of runs in $O$, the second-to-last run of $E$ pairs with the endpoint $\pi_n$ (i.e. creates a chain absorbing it), yielding the additional plus 1 showing up in our bound to account for the final run of $E$. (Note that the number of runs of $O$ cannot exceed that of $E$ due to the comparative number of elements and the alternating behavior.) Otherwise, the last run of $E$ pairs with $\pi_n$.
\end{enumerate}

Counting the partitions formed by this stitching strategy, and adding the necessary offset for endpoints and surplus runs as noted, yields the bound:
\[ P_{max}(n) \le \lfloor\frac{\lceil\frac{n-2}{2}\rceil}{2}\rfloor + 2 \]
This completes the proof.
\end{proof}

\vspace{.5cm}
To illustrate what is happening in the proof, consider the following example:

\begin{example}
Consider the roller coaster permutation $\pi = \{5,3,7,1,8,2,6,4\}$ of length $n=8$.
Here, the endpoints are $\pi_1=5$ and $\pi_8=4$. The internal entries separate strictly into two sets:
\begin{itemize}
    \item \textbf{Bottom Set ($E$):} $\{3, 1, 2\}$ at indices $\{2, 4, 6\}$.
    \item \textbf{Top Set ($O$):} $\{7, 8, 6\}$ at indices $\{3, 5, 7\}$.
\end{itemize}
Decomposing these into increasing runs:
\begin{itemize}
    \item The Bottom ($E$) forms runs: $R_E^{(1)} = \{3\}$ and $R_E^{(2)} = \{1, 2\}$.
    \item The Top ($O$) forms runs: $R_O^{(1)} = \{7, 8\}$ and $R_O^{(2)} = \{6\}$.
\end{itemize}
Applying the stitching strategy from Theorem 3.1:
\begin{enumerate}
    \item The Start ($5$) pairs with the first Top run ($7,8$) $\rightarrow \{5, 7, 8\}$.
    \item The first Bottom run ($3$) pairs with the second Top run ($6$) $\rightarrow \{3, 6\}$.
    \item The second Bottom run ($1, 2$) pairs with the End ($4$) $\rightarrow \{1, 2, 4\}$.
\end{enumerate}
This yields a partition size of 3, which matches the theoretical bound:
$\lfloor\frac{\lceil (8-2)/2 \rceil}{2}\rfloor + 2 = 3$.

\begin{center}
\begin{tikzpicture}[scale=1.2]
    \tikzstyle{vertex}=[circle, draw=black, fill=white, inner sep=2pt]
    \tikzstyle{permline}=[gray!50, thin]
    \tikzstyle{part1}=[blue, thick, ->, >=stealth]
    \tikzstyle{part2}=[red, thick, ->, >=stealth]
    \tikzstyle{part3}=[green!60!black, thick, ->, >=stealth]

    \draw[help lines, color=gray!10] (1,1) grid (8,8);
    
    \coordinate (p1) at (1,5);
    \coordinate (p2) at (2,3);
    \coordinate (p3) at (3,7);
    \coordinate (p4) at (4,1);
    \coordinate (p5) at (5,8);
    \coordinate (p6) at (6,2);
    \coordinate (p7) at (7,6);
    \coordinate (p8) at (8,4);

    \draw[permline] (p1)--(p2)--(p3)--(p4)--(p5)--(p6)--(p7)--(p8);

    \foreach \i in {1,...,8} {
        \node[vertex] (n\i) at (p\i) {};
        \node[below, font=\tiny, gray] at (p\i |- 0,0.8) {\i}; 
    }
    \node[above] at (p1) {5};
    \node[below] at (p2) {3};
    \node[above] at (p3) {7};
    \node[below] at (p4) {1};
    \node[above] at (p5) {8};
    \node[below] at (p6) {2};
    \node[above] at (p7) {6};
    \node[below] at (p8) {4};

    
    \draw[part1] (p1) -- (p3);
    \draw[part1] (p3) -- (p5);
    \node[blue, right] at (8.2, 7) {Part 1: $\{5, 7, 8\}$};

    \draw[part2] (p2) -- (p7);
    \node[red, right] at (8.2, 6) {Part 2: $\{3, 6\}$};

    \draw[part3] (p4) -- (p6);
    \draw[part3] (p6) -- (p8);
    \node[green!60!black, right] at (8.2, 5) {Part 3: $\{1, 2, 4\}$};

    \draw[dashed, gray] (0.5, 4.5) -- (8.5, 4.5);
    \node[gray, font=\small, anchor=east] at (0.5, 7) {\textbf{Top Set ($O$)}};
    \node[gray, font=\small, anchor=east] at (0.5, 2) {\textbf{Bottom Set ($E$)}};

\end{tikzpicture}
\end{center}
\end{example}

To validate the theoretical upper bound established in Theorem 3.1, we conducted an experimental enumeration of partition numbers for small $n$. Using code developed in the Sage computer algebra system, we calculated the exact values of $P_{max}(n)$ for $n < 15$. The results, presented in Table \ref{tab:comparison}, demonstrate that our bound is nearly sharp. Notably, the bound is tight in several cases, however the theoretical bound exceeds the actual value by exactly 1 frequently.

\begin{table}[h]
\centering
\begin{tabular}{ccc}
\toprule
$n$ & Actual $P_{max}(n)$ & Bound $P_{max}(n)$ \\
\midrule
3 & 2 & 2 \\
4 & 2 & 2 \\
5 & 2 & 3 \\
6 & 3 & 3 \\
7 & 3 & 3 \\
8 & 3 & 3 \\
9 & 3 & 4 \\
10 & 3 & 4 \\
11 & 4 & 4 \\
12 & 4 & 4 \\
13 & 4 & 5 \\
14 & 4 & 5 \\
15 & 4 & 5 \\
16 & 4 & 5 \\
17 & 5 & 6 \\
18 & 5 & 6 \\
\bottomrule
\end{tabular}
\caption{Comparison of Experimental Values and Theoretical Bound}
\label{tab:comparison}
\end{table}

\section{A Logarithmic Lower Bound Conjecture}

The experimental data presented in Table \ref{tab:comparison} reveals a striking pattern in the growth of $P_{max}(n)$. While our theoretical upper bound grows linearly ($n/4$), the actual values exhibit logarithmic growth. specifically, observing the doubling indices:
\begin{itemize}
    \item $P_{max}(3) = 2$
    \item $P_{max}(6) = 3$
    \item $P_{max}(12) = 4$
\end{itemize}
This suggests that the partition number satisfies the recurrence $P(2n) = P(n) + 1$, leading to a lower bound of the form $P(n) \ge \lfloor \log_2(n/3) \rfloor + 2$.

We propose that this logarithmic behavior is a direct consequence of the \textit{Recursively Alternating} property established in \cite{Adamczak2017}.
Consider $\pi \in RC(n)$. To minimize the partition size, one might attempt to cover the "Bottom" set (indices $2, 4, \dots$) with a single monotonic subsequence. However, Theorem \ref{thm:recursive} dictates that the remaining "Top" set (indices $1, 3, \dots$) forms a Roller Coaster permutation of size $\lceil n/2 \rceil$ with \textit{reversed} orientation. This flip in alternating direction (from Alternating to Reverse Alternating) at every recursive depth likely prevents monotonic subsequences from efficiently covering multiple layers of the recursion, forcing a new partition for each layer.

\begin{conjecture}[Logarithmic Lower Bound]
For $n \ge 3$, the partition number of a Roller Coaster permutation is bounded below by:
\[ P_{max}(n) \ge \lfloor \log_2(n) \rfloor \]
\end{conjecture}

This conjecture, if proven, would confirm that while Roller Coaster permutations are "complex" in terms of local extrema, they are "simple" in terms of partitioning (logarithmic vs. the $\sqrt{n}$ complexity of general permutations).

\section{Conclusion}

In this paper, we established a theoretical upper bound for the partition number of roller coaster permutations, denoted $P_{max}(n)$. By exploiting the strict value separation between even and odd indices in $RC_n$, we derived the bound:
\[ P_{max}(n) \le \lfloor\frac{\lceil\frac{n-2}{2}\rceil}{2}\rfloor + 2 \]
Experimental values for $n < 15$ (Table \ref{tab:comparison}) indicate that this bound is tight for small $n$.

Future work may focus on narrowing the gap between the theoretical bound and the experimental values. Additionally, the relationship between $P_{max}(n)$ and specific forbidden subsequences warrants further investigation, particularly using the Integer Linear Programming approaches recently proposed by Netto \cite{Netto2022}.


\begin{thebibliography}{99}

\bibitem{Adamczak2017}
Adamczak, W. (2017).
\newblock A Note on the Structure of Roller Coaster Permutations.
\newblock \emph{Journal of Mathematics Research}, \textbf{9}(3), 75--79.

\bibitem{Ahmed2013}
Ahmed, T., \& Snevily, H. (2013).
\newblock Some Properties of Roller Coaster Permutations.
\newblock \emph{Bulletin of the ICA}, \textbf{68}, 55--60.

\bibitem{Aldous1999}
Aldous, D., \& Diaconis, P. (1999).
\newblock Longest increasing subsequences: from patience sorting to the Baik-Deift-Johansson theorem.
\newblock \emph{Bulletin of the American Mathematical Society}, \textbf{36}(4), 413--432.

\bibitem{Andre1879}
Andr\'e, D. (1879).
\newblock D\'eveloppements de sec x et de tang x.
\newblock \emph{Comptes Rendus de l'Acad\'emie des Sciences, Paris}, \textbf{88}, 965--967.

\bibitem{Biniaz2019}
Biniaz, A., Bose, P., Dikgit, K. C., Lubiw, A., Maheshwari, A., M\"utze, T., \& Smid, M. (2019).
\newblock Roller-coasters: Long sequences without short runs.
\newblock \emph{SIAM Journal on Discrete Mathematics}, \textbf{33}(3), 1077--1090.

\bibitem{Botler2021}
Botler, F., \& Netto, B. R. L. (2021).
\newblock New Bounds on Roller Coaster Permutations.
\newblock \emph{Anais do VI Encontro de Teoria da Computa\c{c}\~ao}, 50--53, SBC.

\bibitem{Egge2004}
Egge, E., \& Mansour, T. (2004).
\newblock 132-avoiding Two-stack Sortable Permutations, Fibonacci Numbers, and Pell Numbers.
\newblock \emph{Discrete Applied Mathematics}, \textbf{143}, 78--83.

\bibitem{Kezdy1996}
Kezdy, A. E., Snevily, H. S., \& Wang, C. (1996).
\newblock Partitioning permutations into increasing and decreasing subsequences.
\newblock \emph{Journal of Combinatorial Theory, Series A}, \textbf{73}, 353--359.

\bibitem{Mansour2003}
Mansour, T. (2003).
\newblock Restricted 132-alternating permutations and Chebyshev polynomials.
\newblock \emph{Annals of Combinatorics}, \textbf{7}, 201--227.

\bibitem{Netto2022}
Netto, B. R. L. (2022).
\newblock \emph{Two Problems in Combinatorics: Roller Coaster Permutations \& the Erd\H{o}s-S\'os Conjecture}.
\newblock Master's Thesis, Universidade Federal do Rio de Janeiro.

\bibitem{Stanley2010}
Stanley, R. P. (2010).
\newblock A survey of alternating permutations.
\newblock \emph{Contemporary Mathematics}, \textbf{531}, 165--196.

\bibitem{Stankova1994}
Stankova, Z. E. (1994).
\newblock Forbidden subsequences.
\newblock \emph{Discrete Mathematics}, \textbf{132}, 291--316.

\bibitem{Wagner1984}
Wagner, K. (1984).
\newblock Monotonic coverings of finite sets.
\newblock \emph{Elektronische Informationsverarbeitung und Kybernetik}, \textbf{20}(12), 633--639.

\end{thebibliography}
\end{document}